\documentclass[12pt,a4paper]{amsart}
\usepackage[latin1]{inputenc}
\usepackage{latexsym}
\newtheorem{theorem}{Theorem}

\newtheorem{remark}{Remark}
\newtheorem{lemma}[theorem]{Lemma}

\author{J. M. Almira* and  A. Romero}
\title{Riemannian Geometry and the Fundamental \\ Theorem of Algebra}

\thanks{*Corresponding author}

\thanks{The second author is partially supported by the Spanish MEC-FEDER Grant MTM2010-18099 and the
Junta de Andaluc\'{\i}a Regional Grant P09-FQM-4496.}

\begin{document}

\maketitle

\markboth{J. M. Almira and A. Romero}{Fundamental theorem of algebra}

\thispagestyle{empty}

\begin{abstract}
If a (non-constant) polynomial has no zero, then a certain
Riemannian metric is constructed on the two dimensional sphere.
Several geometric arguments are then shown to contradict this
fact.
\end{abstract}

\vspace*{9mm}


\vspace*{9mm}

\section{Introduction}
In \cite{almira_romero} the authors proved that the Gauss-Bonnet
theorem implies the fundamental theorem of algebra. In this note
we present several new Riemannian geometry arguments which lead
also to the fundamental theorem of algebra. All the proofs are
based on the following technical result:
\begin{lemma} \label{lema}
 If there exists an irreducible polynomial $p(z)$ of degree $n>1$,
then there exists a Riemannian metric $g$ on the sphere
$\mathbb{S}^2$ such that its Gaussian curvature, $K_g$, vanishes
identically.
\end{lemma}
This result was already proved in \cite{almira_romero} but we
state it here in a completely different approach that, in our
opinion, is more systematic than the previous one.

Clearly, the metric stated in Lemma \ref{lema}, if it exists, it
is a quite strange geometric object. Indeed, the second step in
all the proofs we present here consist of showing that this metric
cannot exist. In other words, we will point out several well-known
geometrical obstructions to the construction of a flat metric on
the sphere $\mathbb{S}^2$.

In Section 2 we prove Lemma \ref{lema}. Section 3 is devoted to
explain the distinct arguments leading to a proof of the fact that
sphere is not flat. Finally, in Section 4 we connect our proof to
the field extension version of the fundamental theorem of algebra.

\section{Proof of Lemma \ref{lema}}
Assume $p(z)$ is an irreducible polynomial of degree $n>1$. This
implies that the quotient $A:=\mathbb{C}[z]\,/\langle\,
p(z)\,\rangle $ is a field. Furthermore, the map
$\tau:\mathbb{C}^n\longrightarrow A$ given by
$$\tau(a_0,\cdots,a_{n-1})=a_0+a_1z+\cdots+a_{n-1}z^{n-1}\,+\,\langle\,
p(z)\,\rangle$$defines an isomorphism of complex vector spaces. In
particular,  $$\beta=\big\{\,\tau(0,\cdots,1^{\text{i-th
position}},\cdots,0)\,\big\}_{i=1}^n$$is a basis of $A$. Moreover,
we have that $\tau(-w,1,0,\cdots,0), \tau(-1,w,0,\cdots,0)\neq 0$,
for any $w\in\mathbb{C}$. Hence
$$H(w)=\tau(-w,1,0,\cdots,0)\,\tau(-1,w,0,\cdots,0)\neq 0.$$ Let
$M(w)$ be the associated matrix, with respect to the basis $\beta$
above, to the complex linear operator $L_w:A\longrightarrow A$
given by
$$L_w\big(\tau(a_0,\cdots,a_{n-1})\big)=H(w)\,\tau(a_0,\cdots,a_{n-1}).$$
Obviously, $L_w$ is an isomorphism since $H(w)\neq 0$ and $A$ is a
field. Hence $\det\big(M(w)\big)\neq 0$ for all $w$. Furthermore,
$f(w):=\det\big(M(w)\big)$ is a polynomial.

Now, the linearity of $\tau$ guarantees that, for all
$w\in\mathbb{C}\setminus\{0\}$,
\begin{eqnarray*}
H(1/w) &=& \tau(-1/w,1,0,\cdots,0)\,\tau(-1,1/w,0,\cdots,0)\\[1mm]
&=&
\big[(1/w)\,\tau(-1,w,0,\cdots,0)\big]\,\big[(1/w)\,\tau(-w,1,0,\cdots,0)\big]\\[1mm]
&=& (1/w^2)H(w),
\end{eqnarray*}
so that
$$f(1/w)=\det\big(M(1/w)\big)=\det\big((1/w^{2})M(w)\big)=(1/w^{2n})\det\big(M(w)\big)=(1/w^{2n})f(w).$$

It follows that there exists a Riemannian metric $g$ on the sphere
$\mathbb{S}^2=\widehat{\mathbb{C}}=\mathbb{C}\cup\{\infty\}$, such
that
\[
g=\frac{1}{|f(w)|^{\frac{2}{n}}}\,|dw|^2 \quad \text{for} \quad
w\in \mathbb{C}\]
 and
\[
g=\frac{1}{|f(1/w)|^{\frac{2}{n}}}\,|d(1/w)|^2 \quad \text{for}
\quad w\in \widehat{\mathbb{C}} \setminus\{0\}.
\]

Now, a simple computation shows that the Gaussian curvature $K_g$
of $g$ satisfies
\[
\frac{1}{|f(w)|^{\frac{1}{n}}}\,K_g=\frac{1}{n}\Delta \log
|f(w)|=\frac{1}{n}\Delta \,\mathbf{Re}\log f(w)=0 \quad \text{for
all} \quad w\in\mathbb{C}\setminus \{0\},
\]
since the real part of a holomorphic function must be harmonic.
This obviously implies that $K_g=0$ on the whole sphere and ends
the proof.

\section{The sphere is not flat}

Of course, the following arguments leading the title of this
section are well-known. We recall them for sake of completeness of
this note.

\vspace{2mm}

\noindent {\bf First argument.} Any Riemannian metric on
$\mathbb{S}^2$ must be geodesically complete, from the Hopf-Rinow
theorem. Therefore, the flat Riemannian manifold
$(\mathbb{S}^2,g)$ is also geodesically complete, and, taking into
account that $\mathbb{S}^2$ is connected and simply connected, the
Cartan theorem on the classification of space forms (see
\cite[Theorem 7.10]{KN}, for instance) gives that
$(\mathbb{S}^2,g)$ should be globally isometric to Euclidean
plane, which is impossible because of the compactness of the
sphere.

\vspace{2mm}

\noindent {\bf Second argument.} The usual Riemannian metric of
Gaussian curvature 1 on $\mathbb{S}^2$ is locally written
$g^0=\big(4/(1+|w|^2)^2\big)\, |dw|^2$, $w \in \mathbb{C}$.
Therefore, the Riemannian metric $g$ in Lemma \ref{lema} is
pointwise conformally related to $g^0$, i.e., $g=e^{2u}g^0$, where
$u\in C^{\infty}(\mathbb{S}^2)$ is non constant (note that a
homothetical metric to $g^0$ has constant positive Gaussian
curvature). Using now the relation between the Gaussian curvatures
of two pointwise conformally related metrics we get $\Delta^0 u
=1$, where $\Delta^0$ is the Laplacian relative to the metric
$g^0$. Making use again of the compactness of the sphere, the
classical maximum principle gets that $u$ must be constant which
is impossible.

\vspace{2mm}

\noindent {\bf Third argument \cite{almira_romero}.} As an easy
consequence of the Gauss-Bonnet theorem, any Riemannian metric on
$\mathbb{S}^2$ has some elliptic point, i.e., a point where its Gaussian 
curvature is strictly positive. Hence, the existence of the metric $g$ in Lemma
\ref{lema} contradicts the Gauss-Bonnet theorem.

\begin{remark} {\rm It should be noted that is crucial for our purposes
that the the Gaussian curvature of the metric $g$ in Lemma
\ref{lema} is zero on all $\mathbb{S}^2$. Riemannian metrics on a
sphere with non-constant Gaussian curvature $K$ such that $0 \leq
K \leq 1$ and $K=0$ on a non zero measure set are known to exist.}
\end{remark}

\section{A final comment}
The proof of Lemma \ref{lema} we have presented in this note can
be adapted with no extra effort to give a proof of the fact:
\begin{quote}
If $A$ is a commutative $\mathbb{C}$-algebra, $M$ is a maximal
ideal of $A$ and  $x+M\in A/M$ is algebraic of degree $n>1$ over
$\mathbb{C}$, then there exists a Riemannian metric $g$ on the
sphere $\mathbb{S}^2$ such that $K_g$ vanishes identically.
\end{quote}

This, in conjunction with the arguments in previous section, leads
to a new and direct proof of the following well-known result:
\begin{theorem}[Field extension version of FTA] \label{teo}
Let $A$ be a commutative $\mathbb{C}$-algebra and let $M$ be a
maximal ideal of $A$. If $A/M$ is an algebraic field extension of
$\mathbb{C}$ (in particular, if
$[A/M:\mathbb{C}]=\dim_{\mathbb{C}}(A/M)=n<\infty$, where
$\dim_{\mathbb{C}}V$ denotes complex dimension) then
$[A/M:\mathbb{C}]=1$.
\end{theorem}

\vspace{9mm}

\footnotesize{J. M. Almira

Departamento de Matemáticas. Universidad de Ja\'{e}n.

E.U.P. Linares C/Alfonso X el Sabio, 28

23700 Linares (Ja\'{e}n) Spain

email: {\ttfamily jmalmira@ujaen.es}}

\vspace{3mm}

\footnotesize{A. Romero

Departamento de Geometría y Topología.

Universidad de Granada.

18071 Granada. Spain

email: {\ttfamily aromero@ugr.es}}


\begin{thebibliography}{99}
\bibitem{almira_romero} J.M. Almira, A. Romero, Yet another
application of the Gauss-Bonnet Theorem for the sphere,
\textit{Bull. Belg. Math. Soc. Simon Stevin,} \textbf{14} (2007),
341--342.

\bibitem{KN} S. Kobayashi, K. Nomizu, \textit{Foundations of Differential
Geometry I}, Interscience, New York, 1963.

\end{thebibliography}
\end{document}